\documentclass[preprint,10pt,nopreprintline]{elsarticle}

\usepackage[english]{babel}
\usepackage{geometry}
\usepackage{amsmath,amsfonts,amssymb,amsthm}
\usepackage[final]{hyperref}
\usepackage{graphicx}
\usepackage{empheq}
\usepackage{mathtools}
\usepackage{physics}
\usepackage{mathtools}
\usepackage[dvipsnames]{xcolor}
\usepackage{pdflscape}
\usepackage{etoolbox}
\usepackage{fancyhdr}
\usepackage{lastpage}
\usepackage{ifdraft}
\usepackage{hyphenat}
\usepackage{natbib}
\ifdraft{\usepackage{showlabels}}{}
\usepackage{enumitem}
\usepackage{siunitx}
\usepackage{booktabs}
\usepackage{multirow}
\usepackage{array}
\usepackage{caption}
\usepackage{float}
\usepackage{rotating}
\usepackage[noabbrev]{cleveref}
\usepackage{algorithm}
\usepackage[noend]{algpseudocode}
\newtheorem{example}{Example}

\usepackage{lineno}

\begin{document}

\begin{frontmatter}



\title{Sparse identification of delay equations with distributed memory} 


\author[1]{Dimitri Breda}
\ead{dimitri.breda@uniud.it}
\author[1]{Muhammad Tanveer}
\ead{tanveer.muhammad@spes.uniud.it}
\author[2]{Jianhong Wu}
\ead{wujh@yorku.ca}

\affiliation[1]{organization={CDLab – Computational Dynamics Laboratory Department of Mathematics, Computer Science and Physics – University of Udine},
                addressline={via delle Scienze 206}, 
                city={Udine},
                postcode={33100}, 
                country={Italy}}
\affiliation[2]{organization={Laboratory for Industrial and Applied Mathematics Department of Mathematics and Statistics -- York University},
                addressline={Toronto, ON}, 
                city={Toronto},
                postcode={M3J 1P3}, 
                country={Canada}}

\begin{abstract}
We present a novel extension of the SINDy framework to delay differential equations with {\it distributed delays} and {\it renewal equations}, where typically the dependence from the past manifests via integrals in which the history is weighted through specific functions that are in general nonautonomous. Using sparse regression following the application of suitable quadrature formulas, the proposed methodology aims at directly reconstructing these kernel functions, thereby capturing the dynamics of the underlying infinite-dimensional systems. Numerical experiments confirm the effectiveness of the presented approach in identifying accurate and interpretable models, thus advancing data-driven discovery towards systems with distributed memory.
\end{abstract}




\end{frontmatter}



\section{Introduction} \label{s_{i}ntroduction}
Data-driven discovery of governing equations has become an important approach for studying dynamical systems, revealing their underlying laws based on observational data. A key framework in this area is the Sparse Identification of Nonlinear Dynamics (SINDy)~\cite{bk19,bpk16}. Originally introduced for Ordinary Differential Equations (ODEs), SINDy offers a method for identifying parsimonious, interpretable models from time series data by using sparse regression on a possibly rich library of candidate functions.

The success of SINDy has led to many extensions, including partial differential equations (PDEs)~\citep{demo2018pydmd,rud17}, stochastic differential equations~\citep{boninsegna2018sparse,wanner2024higher}, differential equations with discrete time delays (DDEs)~\citep{bbt24,kopeczi23,pec24,sandoz2023sindy} and, more recently, stochastic delay differential equations \citep{ breda2025sparse,han2024approximation}. These developments demonstrate the flexibility of such an identification framework, solidifying its importance as a key technique in the consolidated, yet growing, field of data-driven paradigms.

However, a significant gap remains in applying SINDy to systems with distributed memory, which more realistically capture past dependence. These dynamics, modeled by delay integro-differential equations (DIDEs) and renewal equations (REs), are central to many biological, ecological, and engineering systems, in which evolution relies on a range of past states~\citep{ando2020fast,beretta2016discrete,bdgpv12,metz2014dynamics}. Distributed delays and renewal-type Volterra integral equations model this memory by integrating over a continuum of past states via suitable kernel functions. This paper represents a first attempt to recover the right-hand side of DIDEs and REs from data by directly identifying these constitutive integral kernels to obtain interpretable models rather than black-box tools to just progress simulation in time. Note that obtaining information on the form of constitutive kernels is essential to advance knowledge in the complex modeling of emerging or strategic fields such as, e.g., behavioural epidemiology \citep{bdm13,buonomo2025minimal,buonomo2025integral}, tick-borne disease transmission \citep{btwz2025tick,CSIAM-LS-1-2} or supply chains for sustainable innovation \cite{guerrini2024dynamics,sarfraz2025leveraging,sipahi2008supply,zhan2022quality}.

Our method is inspired by recent improvements of the SINDy proposed to address the integral formulation of ODEs~\citep{messenger2021weak, wei2022sparse}, introduced mainly for better noise resistance by replacing numerical differentiation with integration. Integrals are instead inherently present in DIDEs and REs, with memory kernels that are typically nonautonomous: they both depend on the delay or integration variable (which might represent, e.g., age in the case of age-structured population models) as well as on system state. This continuous nature of memory integration requires careful handling of numerical quadrature, and the main innovation of our work is a consistent approach for creating libraries of candidate functions that clearly include both the state and delay variables. The proposed framework also adopts improved optimization techniques, such as particle swarm optimization (PS)~\citep{bonyadi2017,kennedy1995,shi1998}, to estimate other essential system parameters such as the relevant integration extrema (e.g., unknown maturation or maximum age of individuals in a population). The resulting tool is thus able to 
simultaneously reveal both the functional form of the memory kernel as well as the integration span.

The contents of the paper are organized as follows. In Section~\ref{sec:dde_re} we introduce the basic necessary mathematical description of DIDEs and REs. In Section~\ref{sec:isindy} we revisit the integral SINDy framework of \citep{bk19,bpk16} for ODEs as a starting basis for the proposed extension. The main methodology for identifying distributed delay systems is introduced in Section~\ref{sec:sindy_kernel}, where we detail the construction of a library based on quadrature and the formulation of sparse regression. Relevant pseudocodes (provided in \ref{sec:appendix}, with MATLAB demo codes freely available at \url{https://cdlab.uniud.it/software}) are also discussed therein, together with an explanatory comparison between the newly proposed approach and a standard ``black-box'' application of SINDy. In Section~\ref{sec:exp} we present extensive numerical experiments to demonstrate the performance of the resulting tool and examine the effects of key algorithm choices. Finally, we briefly comment on the wider implications of this study and outline further future research directions for data-driven discoveries in the field of infinite-dimensional dynamical systems.

\section{Distributed delay differential and renewal equations} \label{sec:dde_re}
In many delayed dynamical systems, the time evolution of the concerned quantity depends not only on its values on just one or few single distinct points in the past, but rather on a weighted contribution of its history (full or possibly partial in between given or unknown time instants), with such weighting described by a delay kernel. These systems are portrayed by models with distributed memory, in which past dependence is expressed via integral functionals. DIDEs and REs are primary classes of such models, together with relevant or even more complex formulations via PDEs (whose treatment we deliberately omit in this work, yet we refer the interested reader to \cite{Boldin02082024,franco2023modelling} as far as models of structured populations are concerned).

In principle, both DIDEs and REs can make use of a common structure for the right-hand side to capture the influence of the history of the system. In this respect a general DIDE is formulated as a differential equation
\begin{equation}\label{eq:dide}
    x'(t) = G(x_t)
\end{equation}
whose rate of change depends on the past, whereas a typical RE
\begin{equation}\label{eq:re}
    x(t) = G(x_t)
\end{equation}
directly relates the current value of the unknown to its history $x_t$ defined as $x_t(\sigma)\coloneqq x(t+\sigma)$, $\sigma\in [-\tau, 0]$, for a given maximum delay $\tau>0$. In both formulations, the central component is a functional $G:X\to\mathbb{R}^{n}$ for $n$ a positive integer, which processes the state history and operates in a Banach space $X$ of functions $[-\tau,0]\to\mathbb{R}^{n}$ of suitable regularity (typically continuous for DIDEs and measurable for REs \cite{dgg07}). $G$ accounts for the aggregated influence at the current time of the history of the system on the time span $[-\tau,0]$\footnote{We omit to consider formulations with multiple integrals as the following treatment can be generalized straightforwardly also thanks to the additivity of integration. Instead in the case of DIDEs we feel free to include next also current time and discrete delay terms.\label{foot1}}. This functional is in general defined by an integral
\begin{equation}\label{eq:G}
    G(\varphi)\coloneqq \int_{-\tau}^{0} g(\sigma, \varphi(\sigma))\dd\sigma,\quad\varphi\in X,
\end{equation}
over the delay interval, where $g:[-\tau,0]\times\mathbb{R}^{n}\rightarrow \mathbb{R}^{n}$ is a smooth possibly nonlinear map termed {\it kernel} function throughout the paper. $g$ determines the specific nature of the model by assigning a weight to past states. The integration variable $\sigma$ represents the history time from the present moment $t$: it directly affects the weight selection (first argument of $g$) beyond the specific moment at which the past is considered (via the second argument of $g$). As we clarify in the following sections, this nonautonomous dependence of $g$ from $\sigma$ is a key aspect for its reconstruction from observed data. We instead conclude this section with a couple of prototype examples of DIDEs and REs.


\begin{example}
The logistic Daphnia model \citep{breda2016pseudospectral,breda2020approximation} 
\begin{subequations}\label{eq:daphnia_re_dd} 
\begin{empheq}[left=\empheqlbrace]{align}
b(t) &= \beta S(t)\!\int_{a_*}^{a_\dagger} b(t-a)\,da, 
\label{eq:daphnia_re_dd_a} \\[1ex]
S'(t) &= r S(t)\!\left(1-\frac{S(t)}{\mathcal{K}}\right)
        - \gamma S(t)\!\int_{a_*}^{a_\dagger} b(t-a)\,da,
\label{eq:daphnia_re_dd_b}
\end{empheq}
\end{subequations}
offers a simplified framework of an RE coupled to a DIDE to describe a consumer-resource dynamics. Above $b$ denotes the birth rate of a population feeding on a resource $S$ that in the absence of the consumer grows logistically with the growth rate $r$ and carrying capacity $\mathcal{K}$. $\beta$ and $\gamma$ are the effective fertility and consumption rates of the adult individuals of the consumer population, which are assumed to have a constant survival probability between the maturation age $a_{\ast}$ and the maximal age $a_{\dag}$. In Section \ref{sec:exp} we employ \eqref{eq:daphnia_re_dd} as a representative example of models of coupled DIDEs/REs, making it an ideal starting point for demonstrating the effectiveness of the proposed sparse identification of distributed memory kernels. In the experiments, we will treat both the maturation and the maximum ages as unknown. In terms of \eqref{eq:dide} and \eqref{eq:re} (and recall footnote \footref{foot1}) the right-hand side of \eqref{eq:daphnia_re_dd} can be rewritten as
\begin{equation*}
f(b_{t},S_{t})+G(b_{t},S_{t})
\end{equation*}
with
\begin{equation*}
f(\varphi,\psi)\coloneqq\left(0,r\psi(0)\left(1-\dfrac{\psi(0)}{\mathcal{K}}\right)\right)^{T},\qquad
G(\varphi,\psi)\coloneqq\int_{-a_{\dag}}^{-a_{\ast}}g(\sigma,\varphi(\sigma),\psi(\sigma))\dd\sigma
\end{equation*}
for
\begin{equation*}
g(\sigma,\varphi(\sigma),\psi(\sigma))\coloneqq\left(\beta\psi(0)\varphi(\sigma),-\gamma\psi(0)\varphi(\sigma)\right)^{T}.
\end{equation*}
In Section \ref{sec:exp} we will identify both $f$ and $g$. Note that in this case, the kernel $g$ is autonomous.
\end{example}
\begin{example}
In the behavioural epidemiology framework considered, e.g., in \citep{bdm13}, the classical SIR paradigm can be expanded to include vaccination conditioned on information decision processes. Consider a stable and equally interacting population structured into compartments comprising susceptible ($S$), infected ($I$), recovered ($R$) and vaccinated ($V$) individuals, with the vaccination rate governed by an information variable $M(t)$ tracking past disease levels. A plausible mathematical model reads
\begin{equation*}
\left\{\setlength\arraycolsep{0.1em}\begin{array}{rcl}
S'(t) &=& \mu \left[1 - p(M(t))\right] - \mu S(t) - \beta S(t) I(t), \\
    I'(t) &=& \beta S(t) I(t) - (\mu+\nu) I(t), \\
    R'(t) &=& \nu I(t) - \mu R(t), \\
    V'(t) &=& \mu p(M(t)) - \mu V(t),
\end{array}\right.
\end{equation*}
where $\mu$ is the demographic turnover rate, $\beta$ is the contact rate, and $\nu$ is the recovery rate. The coverage function $p(M)$ is assumed to be non-decreasing and is used to express the proportion of immunized newborns in terms of an information variable
\begin{equation*}
    M(t)\coloneqq \int_{-\infty}^{0} K(\sigma)h(S(t+\sigma), I(t+\sigma))\dd\sigma.
\end{equation*}
Here, $h:[0,+\infty)\times [0,+\infty) \rightarrow \mathbb{R}$ represents the prevalence of the disease weighted through a normalized memory kernel $K$. Typically $K$ is an Erlang distribution, but realistic data may require other choices asking for suitable identification tools. Note that in general the local stability at the endemic equilibrium is affected by the properties of $K$ \cite{ando2020fast}: as the memory concentration increases, the system is likely to have an increased number of intervals with instabilities like Hopf bifurcations and oscillations. This model involves an infinite delay and although analytical tools can handle this case in particular situations (e.g., using the linear chain trick \cite{mac78} as in \cite{dms07}) as well as suitable quadrature rules on the half-line can be employed, time truncation can proves effective for identification, as we will show in Section \ref{sec:exp} for another (Ricker-type) population model with infinite delay that we will consider first in Section \ref{sec:ricker}.
\end{example}
\section{Integral SINDy for ODEs}\label{sec:isindy}
The SINDy framework of \citep{bk19,bpk16} for ODEs recovers parsimonious governing equations from time-series data for dynamical systems. When time derivatives are unavailable, the standard approach relies on numerical differentiation, which may lead to significant sensitivity to the measurement noise. The alternative application of SINDy to the integral form of the underlying ODE has been proposed to improve robustness \citep{messenger2021weak,schaeffer2017sparse,wei2022sparse}. In this version, hereafter referred to as integral SINDy, the derivative estimation is replaced with numerical integration, which is an inherently more stable operation in noisy settings. This choice has inspired our treatment of systems with distributed memory, therefore we briefly describe its main features next.

Consider the IVP
\begin{equation}\label{eq:ode_system}
\left\{\setlength\arraycolsep{0.1em}\begin{array}{rcl}
x'(t) &=& f(x(t),\quad t\in[0,T],\\
x(0)&=&x_{0},
\end{array}\right.
\end{equation}
for a general autonomous system of ODEs where ${f}: \mathbb{R}^{n} \to \mathbb{R}^{n}$ is an unknown vector field, $x_{0}\in\mathbb{R}^{n}$ and $T>0$. Integrating ~\eqref{eq:ode_system} yields the equivalent integral form
\begin{equation}\label{eq:integral_form}
x(t) - x_0 = \int_{0}^{t} {f}\left({x}(\sigma)\right) \dd\sigma,
\end{equation}
In practical scenarios, the state is observed at discrete times $t_{1},\ldots,t_{m}$ and the corresponding snapshots are organized in the data matrix
\begin{equation}\label{eq:data}
    \mathbf{X} =(X_{1},\ldots,X_{n})
\in \mathbb{R}^{m \times n}.
\end{equation}
Each component $f_{j}$ of the vector field is approximated as a sparse linear combination of candidate functions from a predefined library ${\Theta}(x)\in\mathbb{R}^{m\times p}$. The library may include polynomial, trigonometric, or other nonlinear terms. Hence
\begin{equation}
f_{j}({x}) \approx {\Theta}(x^{T}) {\xi}_{j},
\end{equation}
where ${\xi}_{j}\in\mathbb{R}^{p}$ is a sparse coefficient vector for the $j$-th governing equation.
Substituting the latter into \eqref{eq:integral_form} requires integrating the candidate functions of the library. To this aim we resort to a simple rectangles formula based on the same time instants of the data sampling, assumed to be uniformly distributed in $[0,T]$ for the sake of simplicity. By letting $h\coloneqq t_{k}-t_{k-1}$ we get
\begin{equation}\label{eq:discrete_form}
 x_{j}(t_i) - x_{0,j} \approx \left(\int_{0}^{t_i} {\Theta}(x^{T}(\sigma))\dd\sigma\right){\xi}_{j} \approx \left(h\sum_{k=0}^{i-1} {\Theta}(x^{T}(\sigma_{k}))\right) {\xi}_{j}
\end{equation}
for each state component $j=1,\ldots,n$ and every time instant $t_{i}$, $i=1,\ldots,m$.
This establishes a relationship between the change in the state variable over the interval $[0,T]$ and a linear combination of the functions of $\Theta$ through the unknown coefficients ${\xi}_{j}$. According to the standard SINDy procedure, the regularized least squares regression
\begin{equation}\label{eq:ISINDYreg}
\hat{{\xi}}_{j} = \arg\min_{\xi\in\mathbb{R}^{p}}\left( \left\|\mathbf{X}_{j}-x_{0,j}-\left(h\sum_{k=0}^{i-1}  {\Theta}(\mathbf{X}_{k,:})\right)\xi\right\|_{2}+\lambda\|\xi\|_{1}\right)
\end{equation}
is solved independently for each $j=1,\ldots,n$ to find the sparse vector ${\xi}_{j}$ that best describes the dynamics $f_{j}$. Above $\mathbf{X}_{k,:}$ denotes the $k$-th row of $\mathbf{X}$ and $\lambda$ is a regularization parameter that controls sparsity. \eqref{eq:ISINDYreg} is a standard convex optimization problem that can be solved efficiently using established tools such as STLS \citep{bk19} or LASSO \cite{lasso96}. In general one can resort to other quadrature formulas and the approach can be easily adapted to the case of non-equispaced samples by using (piecewise) linear interpolation. Next we adopt the underlying use of quadrature to tackle equations with distributed delays, where in general a remarkable difference from what described here is that the integral kernel is by nature {\it nonautonomous}.
\section{Extending SINDy to identify distributed memory kernels}\label{sec:sindy_kernel}
We first illustrate in Section \ref{sec:quadrature} our general approach to identify the kernel $g$ of a distributed memory dependence that characterises the right-hand side \eqref{eq:G} of a DIDE \eqref{eq:dide} or an RE \eqref{eq:re}. Then in Section \ref{sec:optimization} we briefly discuss relevant algorithms including the use of advanced optimization tools to recover possible unknown delays or non-multiplicative parameters\footnote{Multiplicative ones are indeed contained in the regression vectors $\xi_{j}$.} . In the sequel, we use the term ``external'' optimization to refer to this specific problem, which in this context is mainly devoted to recovering unknown integration extrema beyond model parameters characterising special nonlinearities that might be included in the SINDy library in a physics-informed fashion. Finally, in Section \ref{sec:ricker} we compare the new approach with a standard application of SINDy aimed at just reconstructing the unknown dynamics as a ``black-box'' model, which in general does not provide any information on the underlying integration kernel. This analysis is performed on a specific DIDE with a nonlinear and nonautonomous kernel of Ricker type.

\subsection{Quadrature and weighted libraries}
\label{sec:quadrature}
In order to extend SINDy to recover distributed memory right-hand sides, it is key to first apply a suitable quadrature formula to approximate the integral representing $G(x_{t})$ in \eqref{eq:G} as a sum of weighted kernel values:
\begin{equation}\label{eq:quadrature}
\int_{-\tau}^{0} g(\sigma,x(t+\sigma))\dd\sigma\approx\sum_{k=1}^{K} w_{k} g(\sigma_{k},x(t+\sigma_{k})). 
\end{equation}
The above $K$ represents the number of quadrature nodes $\sigma_{k}\in[-\tau, 0]$ and corresponding weights $w_{k}$, $k=1,\ldots,K$. According to the driving principle behind SINDy, we assume that the kernel $g$ admits a sparse linear representation within a library of candidate functions $\Theta$. Since at each time $t$ the kernel $g$ is a function of both the delay variable $\sigma$ and the state value $x(t+\sigma)$, the library must take into account both arguments, thus allowing the recovery of {\it nonautonomous} kernels. Therefore, the delay variable $\sigma$ is treated as an explicit candidate function, together with its possible nonlinear dependencies. Then, in view of the quadrature \eqref{eq:quadrature}, for each $k=1,\ldots,K$ we first consider a library $\Theta_{k}$ of candidate functions $p$ which depend on the quadrature node $\sigma_{k}$ and time-shifted state values $x(t+\sigma_{k})$\footnote{Note that in principle $p$ can depend on $k$ even though we do not resort to such generality here.}. Following the usual SINDy description relying on the available time samples organized in the matrix $\mathbf{X}$ as in \eqref{eq:data}, such a $k$-th library may typically include polynomial terms up to a given degree $d$ as
\begin{equation*}
\setlength\arraycolsep{0.1em}\begin{array}{rcl}
\Theta_{k}\coloneqq\Theta(\sigma_{k},\mathbf{X}(\cdot+\sigma_{k}))&=\big[&1,\sigma_{k},\mathbf{X}_{1}(\cdot+\sigma_{k}),\ldots,\mathbf{X}_{n}(\cdot+\sigma_{k}),\\[2mm]
&&\sigma_{k}^{2},\sigma_{k}\mathbf{X}_{1}(\cdot+\sigma_{k}),\ldots ,\sigma_{k}\mathbf{X}_{n}(\cdot+\sigma_{k})),\mathbf{X}_{1}^{2}(\cdot+\sigma_{k}),\\[2mm]
&&\ldots,\sigma_{k}^{d},\sigma_{k}^{d-1}\mathbf{X}_{1}(\cdot+\sigma_{k}),\ldots,\mathbf{X}_{n}^{d}(\cdot+\sigma_{k})\big],
\end{array}
\end{equation*}
plus any other nonlinear function of either $\sigma_{k}$ or $\mathbf{X}_{j}(\cdot+\sigma_{k})$, $j=1,\ldots,n$, that the user can add based on prior model knowledge or any other available information (e.g., $e^{\sigma_{k}}$, $\sin(\mathbf{X}_{j}(\cdot+\sigma_{k}))$, etc.). Above $\mathbf{X}_{j}(\cdot+\sigma_{k})$ represents the j-th state variable evaluated at shifted time instants $t_{i}+\sigma_{k}$, $i=1,\ldots,m$: such values can be recovered by (piecewise linear) interpolation if not already at disposal within $\mathbf{X}$. 

A weighted sum of these libraries $\Theta_{k}$, $k=1,\ldots,K$, through the chosen quadrature weights $w_{k}$ represents the complete SINDy library to be used in the final regression problem. For an RE \eqref{eq:re} the corresponding linear representation reads
\begin{equation}\label{eq:reg-re}
\mathbf{X}_{j} \approx\left( \sum_{k=1}^{K} w_{k} \Theta(\sigma_{k}, \mathbf{X}(\cdot+\sigma_{k}))\right)\xi_{j},
\end{equation}
for each $j$-th component of the state, $j=1,\ldots,n$. For a DIDE \eqref{eq:dide} it is enough to replace $\mathbf{X}_{j}$ at the left-hand side with samples data $\mathbf{X}'_{j}$ for the time derivatives\footnote{Derivative values can be recovered from $\mathbf{X}$ via suitable techniques for noise reduction~\cite{chartrand2011numerical} if not available from measurements. Alternatively, one can resort to the integral reformulation of the DIDE as explained in Section \ref{sec:isindy} for ODEs.}. In \eqref{eq:reg-re} the vector $\xi_{j}\in\mathbb{R}^{p}$ is the sparse vector of coefficients sought that finally selects the active library terms for the $j$-th equation. To get each $\xi_{j}$ we independently solve the regularized regression problem
\begin{equation*} 
\hat{\xi}_{j} = \arg\min_{\xi\in\mathbb{R}^{p}} \left(\left\| \mathbf{X}_{j} - \left( \sum_{k=1}^{K} w_{k} \Theta(\sigma_{k}, \mathbf{X}(\cdot+\sigma_{k}))\right)\xi\right\|_{2}+\lambda\|\xi\|_{1}\right) ,
\end{equation*}
for each $j=1,\ldots,n$. By collecting all the resulting sparse vectors in a matrix $\Xi\coloneqq (\hat{\xi}_{1},\ldots,\hat{\xi}_n)\in\mathbb{R}^{p\times n}$, we finally identify \eqref{eq:re} with
\begin{equation*}
x(t)\approx\left(\int_{-\tau}^{0}\Theta(\sigma,x_{t}^{T})\dd\sigma\right)\Xi
\end{equation*}
(or \eqref{eq:dide} with $
x'(t)$ replacing $x(t)$ at the left-hand side), thus revealing the underlying analytical structure of the kernel $g$ in \eqref{eq:G}.
\subsection{Algorithms and external optimization}
\label{sec:optimization}
Here, we provide a brief description of the algorithms included in Appendix \ref{sec:appendix} that are used in Section \ref{sec:exp} for the numerical tests.

Algorithm \ref{alg:distributed-delay-sindy} outlines a systematic method for identifying DIDEs \eqref{eq:dide} and REs \eqref{eq:re} using the proposed SINDy framework as illustrated in the previous section. The resulting tool relies on a priori knowledge of the involved delays, which give indeed the concerned integration window(s). In realistic situations this information is hardly available, and the same consideration holds for any other (non-multiplicative) parameter characterizing possible specific nonlinear dependencies (e.g., the exponent of an exponential function or the frequency or phase of a sinusoid), as well as for numerical parameters such as the number $K$ of quadrature nodes or their distribution. For this reason, following \cite{bbt24}, we integrate the new SINDy framework of Algorithm \ref{alg:distributed-delay-sindy} into a second algorithm, viz. Algorithm~\ref{alg:adaptive-distributed-sindy}, which performs an iterative external optimization loop over the parameters of interest to minimize the SINDy reconstruction error given by Algorithm \ref{alg:distributed-delay-sindy}. We introduce the latter for REs as
\begin{equation}\label{eq:optim_error}
\varepsilon({\rho})\coloneqq\|\mathbf{X}-\Theta(\rho)\Xi(\rho)\|_{2},
\end{equation}
where $\rho$ denotes the vector of possible parameter values in the current step of Algorithm~\ref{alg:adaptive-distributed-sindy} (for DIDEs replace again $\mathbf{X}$ with $\mathbf{X}'$)\footnote{Note that in the experiments part of the data are used to train SINDy and part to validate the recovered model; In this respect \eqref{eq:optim_error} is extended to account for both contributions.}.

External optimization methods, such as PS~\citep{bonyadi2017,kennedy1995,shi1998} or Bayesian optimization (BO)~\citep{garnett2023,shahriari2015,williams2006}, are used to effectively search the parameter space to reduce the error function to a value such that both the recovered kernel and its associated parameters are properly aligned with the observed data. Therefore Algorithm~\ref{alg:adaptive-distributed-sindy} helps uncover not only the sparse structure of the distributed delay kernel (via the ``internal'' regression of SINDy) but also its bounds, other embedded parameters, and possible numerical settings (via the ``external'' optimization loop), thereby enhancing both the interpretability and the predictive accuracy of the identified model.
\subsection{Comparing kernel reconstruction to a black-box SINDy approach}
\label{sec:ricker}
To illustrate the differences between applying the proposed specialized framework for distributed delays (hereinafter named DD-SINDy) and a standard SINDy ``black-box'' approach (hereinafter named BB-SINDy) we consider the DIDE
\begin{equation}\label{eq:ricker}
x'(t) = d_{0}\left(\int_{-\infty}^{0} F_{n/\tau}^{(n)}(\sigma)e^{\eta+d_{1}\sigma-ax(t+\sigma)}x(t+\sigma)\dd\sigma-x(t)\right)
\end{equation}
modeling the growth of a single-species population according to a gamma-distributed maturation delay and a Ricker-type nonlinear reproduction with positive parameters $\eta$ and $a$ (see \cite[equation (4)]{beretta2016discrete} or \cite{jiang2010bifurcation,ricker1954stock}). $x(t)$ denotes adult population size at time $t\geq0$, $d_{0}>0$ the constant death rate of the population, and $d_1\geq0$ the constant death rate during the maturation stage. The integral term represents the history-dependent birth rate of the population, with a delay continuously distributed according to
\begin{equation}\label{eq:gamma}
F_{\alpha}^{(n)}(\sigma)\coloneqq\frac{\alpha^{n}(-\sigma)^{n-1} e^{\alpha\sigma}}{(n-1)!},\quad\sigma\in(-\infty,0],
\end{equation}
for $n$ a positive integer and $\alpha\coloneqq n/\tau$, with mean $\tau>0$ and variance $\tau^{2}/n$. Our objective is to recover the nonautonomous kernel
\begin{equation}\label{eq:rickerg}
g(\sigma,\varphi(\sigma))\coloneqq d_{0}F_{\alpha}^{(n)}(\sigma) e^{\eta+d_{1}\sigma-a\varphi(\sigma)}\varphi(\sigma)
\end{equation}
using the time-series data of $x(t)$. For the sake of focusing on the comparison of the DD-SINDy with BB-SINDy, in this section we deliberately tackle a simple instance of \eqref{eq:ricker} by choosing $d_{1}=0$, $d_{0}=a=1$, $\eta=\log 80$, $n=4$ and $\tau=1$, in which case $\alpha=4$ and the explicit kernel \eqref{eq:rickerg} reduces to
\begin{equation}\label{eq:rickerg1}
g(\sigma,\varphi(\sigma))\coloneqq\gamma\cdot\sigma^{3} e^{4\sigma}e^{-\varphi(\sigma)}\varphi(\sigma)
\end{equation}
for $\gamma\coloneqq -10\,240/3$. The above choices return $g(\sigma,x(t+\sigma))$ as a product of integer powers of the basis functions $\sigma$, $e^{\sigma}$, $e^{-x(t+\sigma)}$, $x(t)$ and $x(t+\sigma)$, which we therefore use as candidates to construct a polynomial SINDy library $\Theta_{{\rm DD}}$ of degree (at least) $4$ for DD-SINDy. Moreover, we further simplify model \eqref{eq:ricker} by truncating the influence of the history on the integration window $[-10\tau,0]$ (in Section \ref{sec:exp} we consider a more general case with non-integer powers and irrational exponents to challenge the external optimization Algorithm~\ref{alg:adaptive-distributed-sindy} in combination with kernel recovery via DD-SINDy). 

As far as BB-SINDy is concerned, \eqref{eq:ricker} is seen just as an instance of the general form \eqref{eq:dide}, so that we seek for a sparse representation of the right-hand side as $G(x_{t})\approx\Theta_{{\rm BB}}(x_{t}^{T})\Xi.$
By resorting to a general polynomial library it is clear that one expects BB-SINDy to correctly recover the coefficient $-d$ of the linear term in \eqref{eq:ricker}, but in general all the other coefficients will not serve to identify interpretable structures for the right-hand side, although the returned model can well serve for making predictions via simulation beyond the time horizon of available data.

\begin{figure}[h!]
\centering
\includegraphics[width=0.65\linewidth]{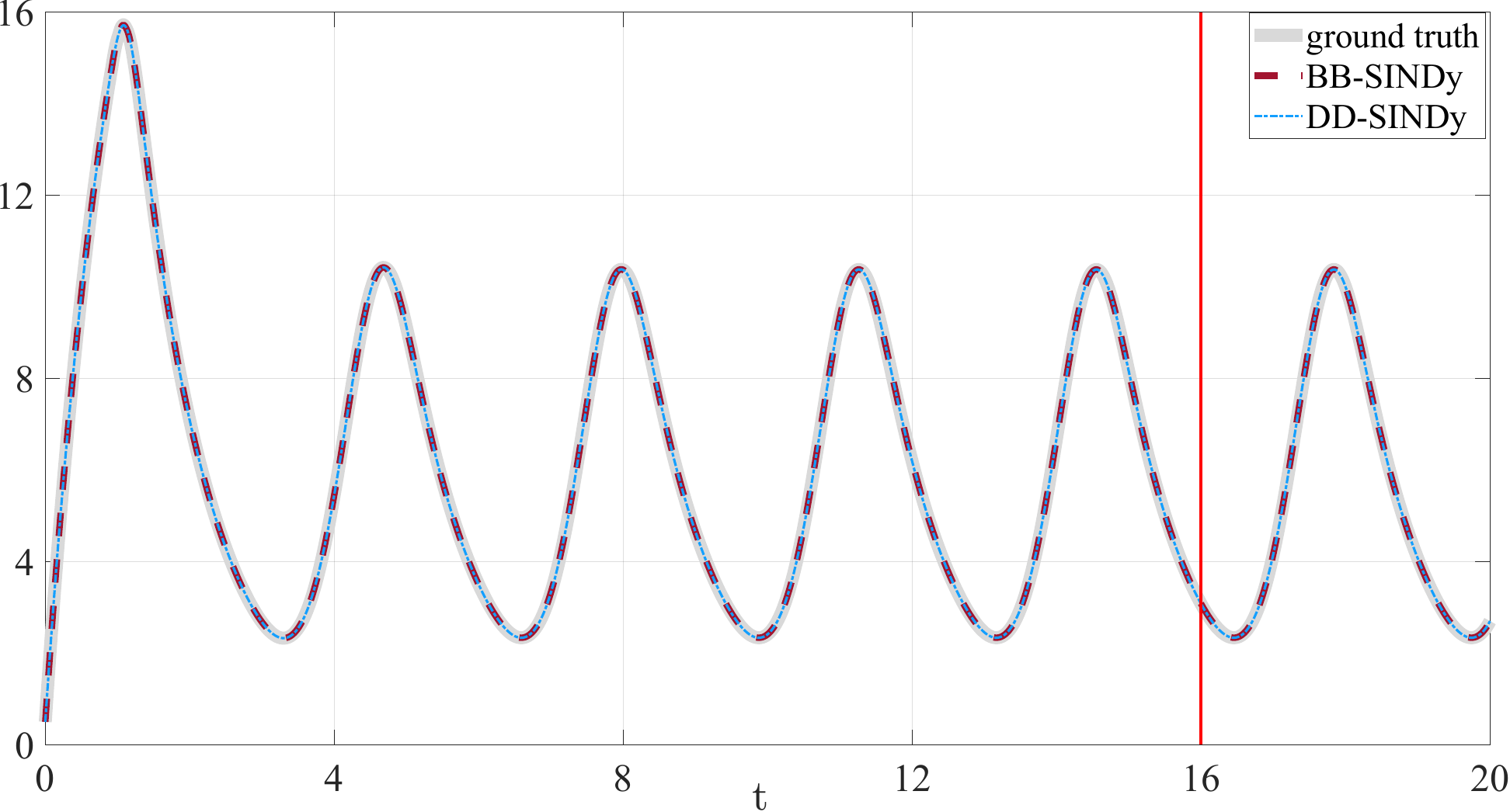}
\caption{comparison of trajectories of \eqref{eq:ricker} reconstructed by DD-SINDy and BB-SINDy. $m=100$ uniform samples were collected in the time window $[0,20]$ with $K=50$ and $\lambda =10^{-2}$, using $80\%$ for training and $20\%$ for validation.}\label{fig:gamma_DD}
\end{figure}

For the sake of comparison synthetic time series data for the population $x$ were generated in the time window $[0,20]$ by numerically solving the IVP for \eqref{eq:ricker} with initial function $\varphi(s)=0.5$, $s\in[-10\tau,0]$. Actually MATLAB's \texttt{dde23} was used to integrate a version of \eqref{eq:ricker} with the integral replaced by a uniform rectangles-based quadrature with $K = 50$ uniform nodes. Then $100$ samples were selected, with $80\%$ reserved for training both DD-SINDy and BB-SINDy, while the remaining data were used for validation to assess the accuracy of the prediction. The same number and distribution of quadrature nodes was used for DD-SINDy, while the sparsity-promoting hyperparameter was set to $\lambda = 10^{-2}$ for both DD-Sindy and BB-SINDy. STLS was used to solve the relevant sparse regression problem. The trajectories of the two recovered models are compared to the ground truth in Figure~\ref{fig:gamma_DD}, showing a remarkable reconstruction capacity for both training and validation.

The root mean square error (RMSE) \eqref{eq:optim_error}measured on the derivatives was also used to evaluate the performance of both approaches, see Table~\ref{tab:comparisondl} . The latter shows also the accuracy of the coefficients learnt during sparse regression. DD-SINDy is able to recover both $d$ in \eqref{eq:ricker} and $\gamma$ in \eqref{eq:rickerg1}, while BB-SINDy returns only the former as already explained.
\begin{table}[htbp]
\centering
\small 
\begin{tabular}{llrrr}
 \toprule
& $|d_{0}-\widehat{d}_{0}|$ & $|\gamma-\widehat{\gamma}|$&  RMSE$_{x^\prime}$ {training} & RMSE$_{x^\prime}$ {validation } \\

 \midrule
 DD-SINDy &   $3.29\times10^{-3}$ & $4.92\times10^{-3}$  & $4.36\times10^{-2}$ & $3.46\times10^{-2}$\\
 BB-SINDy &  $6.81\times10^{-2}$  & -& $3.51\times10^{-2}$ & $9.96\times10^{-2}$\\ 

 
 \bottomrule
 \end{tabular}%
 \caption{comparison of DD-SINDy with BB-SINDy in terms of the absolute error on the reconstructed coefficients (top) and RMSE \eqref{eq:optim_error} on derivatives (bottom) for \eqref{eq:ricker} in the time window $[0, 20]$ using $m=100$ samples and $\lambda=10^{-2}$, with DD-SINDy using $K=50$ quadrature nodes.}
 \label{tab:comparisondl}
 \end{table}
\section{Numerical experiments}\label{sec:exp}
In this section we present the results of quantitatively assessing the proposed methodology on three different models by using Algorithm \ref{alg:distributed-delay-sindy} and Algorithm \ref{alg:adaptive-distributed-sindy} discussed in Section \ref{sec:optimization} and implemented in Appendix \ref{sec:appendix}.
\subsection{A logistic renewal equation}\label{sec:exp_logistic}
\begin{figure}[h!]
\centering
\includegraphics[width=0.65\linewidth]{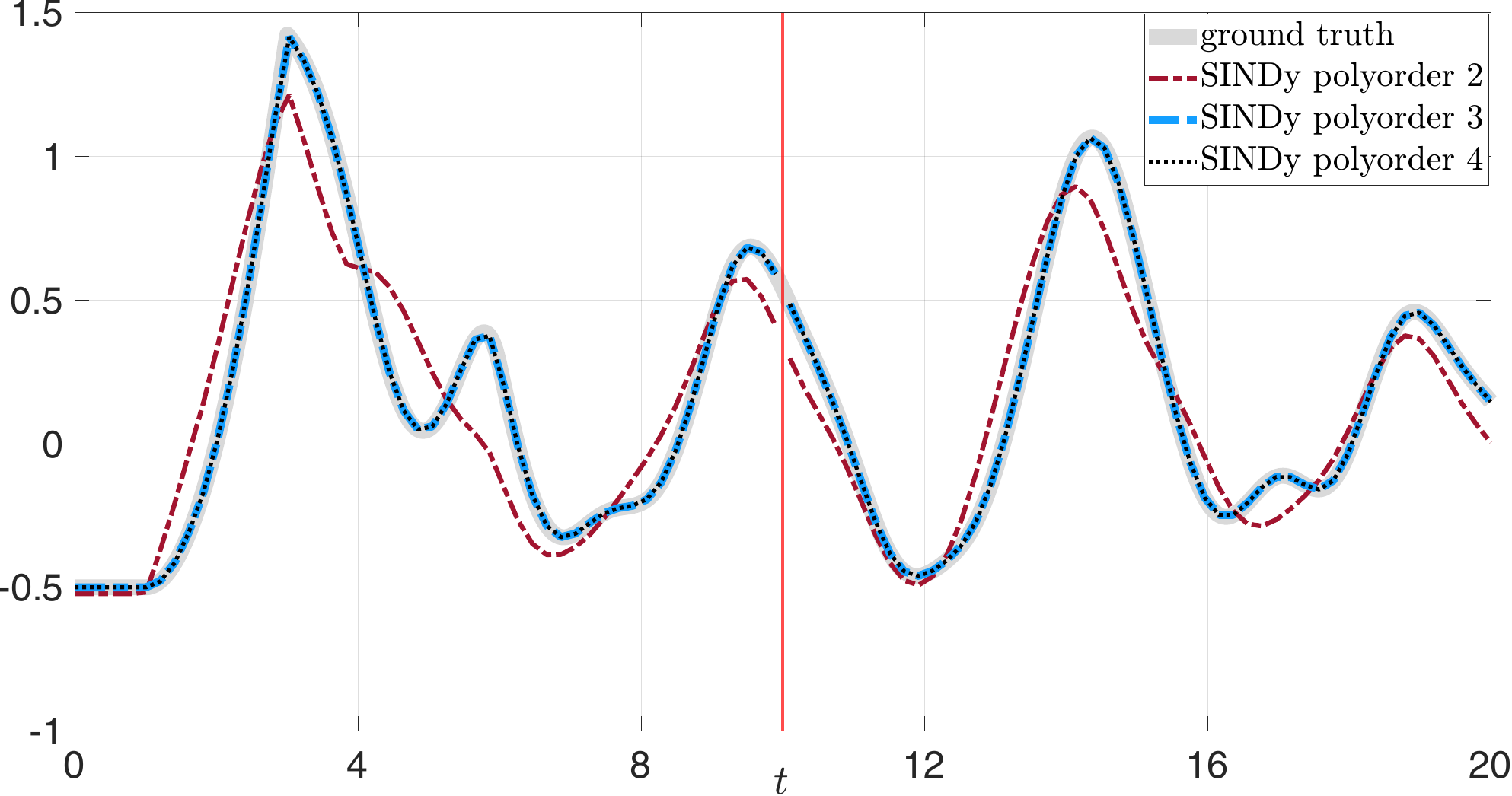}
\caption{comparison of trajectories of \eqref{eq:toy_re} reconstructed by DD-SINDy with a polynomial library of degree $d=2$, $3$ and $4$; corresponding data in Table \ref{tab:combined_results}.}\label{fig:TM1}
\end{figure}

We first consider the prototype RE
\begin{equation}\label{eq:toy_re}
     x(t) = \int_{-3}^{-1} (\sigma+1) x(t+\sigma) [1 - x(t+\sigma)] \,\dd\sigma
\end{equation}
as a toy modification of a classical logistic integral model (see, e.g., \cite{breda2016numerical,breda2018approximation,breda2023practical}.
Data were generated using a time integration scheme from \cite{messina2008stable} on the relevant IVP with constant history function $\varphi(s)= 0.5$, $s\in [-3,0]$. $m$ equidistant samples were collected on $[0, 20]$, using the first $50\%$ for training and the rest for validation. Polynomial SINDy libraries of degree $d=2$, $3$ and $4$ were constructed in combination with trapezoidal quadrature. Table \ref{tab:combined_results} presents the active library terms recovered by DD-SINDy with $m = 100$, $K = 128$ and $\lambda=10^{-5}$, together with relevant RMSE on the derivative for both training and validation. Correspondingly, Figure \ref{fig:TM1} illustrates the reconstructed trajectories. The quadratic library proves insufficient to capture the system dynamics, while extending the library to $d=4$ slightly degrades parameter accuracy due to the inclusion of redundant higher order terms. Conversely, the cubic library achieves the best performance, reducing the coefficient error to at least $10^{-3}$.
\begin{table}[htbp]
\centering
\small
\begin{tabular}{lrrr}
\toprule
polynomial degree $d$
& 2 & 3 & 4 \\
\midrule
$x(t+\sigma)$ &  $6.11\times 10^{-1}$ & $1.82\times 10^{-3}$ &$2.69\times 10^{-1}$  \\
$\sigma\cdot x(t+\sigma)$ &  $6.90\times 10^{-1}$ &$3.06\times 10^{-3}$ & $4.62\times 10^{-1}$  \\
$x(t+\sigma)^{2}$ &  $4.19\times 10^{-2}$ &$2.74\times 10^{-3}$ & $1.30\times 10^{-2}$ \\
$\sigma^{2}$ & $9.89\times 10^{-1}$ &- &  -  \\
$\sigma^{2} \cdot x(t+\sigma)$ &  - & $9.00\times 10^{-3}$ & $2.55\times 10^{-1}$  \\
$\sigma\cdot x(t+\sigma)^{2}$ &  - & $3.55\times 10^{-4}$ & $9.93\times 10^{-3}$  \\
$x(t+\sigma)^3$ &  - & - & $5.13\times 10^{-3}$  \\
$\sigma^3\cdot x(t+\sigma)$ & - & - & $4.19\times 10^{-2}$  \\
$\sigma^{2} \cdot x(t+\sigma)^{2}$ & - & - & $1.64\times 10^{-3}$  \\
$\sigma\cdot x(t+\sigma)^3$ & - & - & $1.90\times 10^{-3}$  \\
$x(t+\sigma)^4$ & - & - & $8.88\times 10^{-4}$ \\
\toprule
RMSE$_{x^\prime}$ training & $2.01\times 10^{-1}$ & $1.40\times 10^{-3}$ & $1.25\times 10^{-3}$ \\
RMSE$_{x^\prime}$ validation &  $1.44\times 10^{-1}$ & $2.18\times 10^{-4}$ & $3.70\times 10^{-4}$ \\
\bottomrule
\end{tabular}
\caption{absolute errors on the coefficients of the polynomial library terms recovered by DD-SINDy on \eqref{eq:toy_re} with degree $d=2$, $3$ and $4$ (top); relevant RMSE on the derivative (bottom). $m=100$ uniform samples were collected in the time window $[0,20]$ with $K=128$ and $\lambda =10^{-5}$, using $50\%$ for training and $50\%$ for validation.}\label{tab:combined_results}
\end{table}

\begin{figure}
    \centering
    \includegraphics[width=1\textwidth]{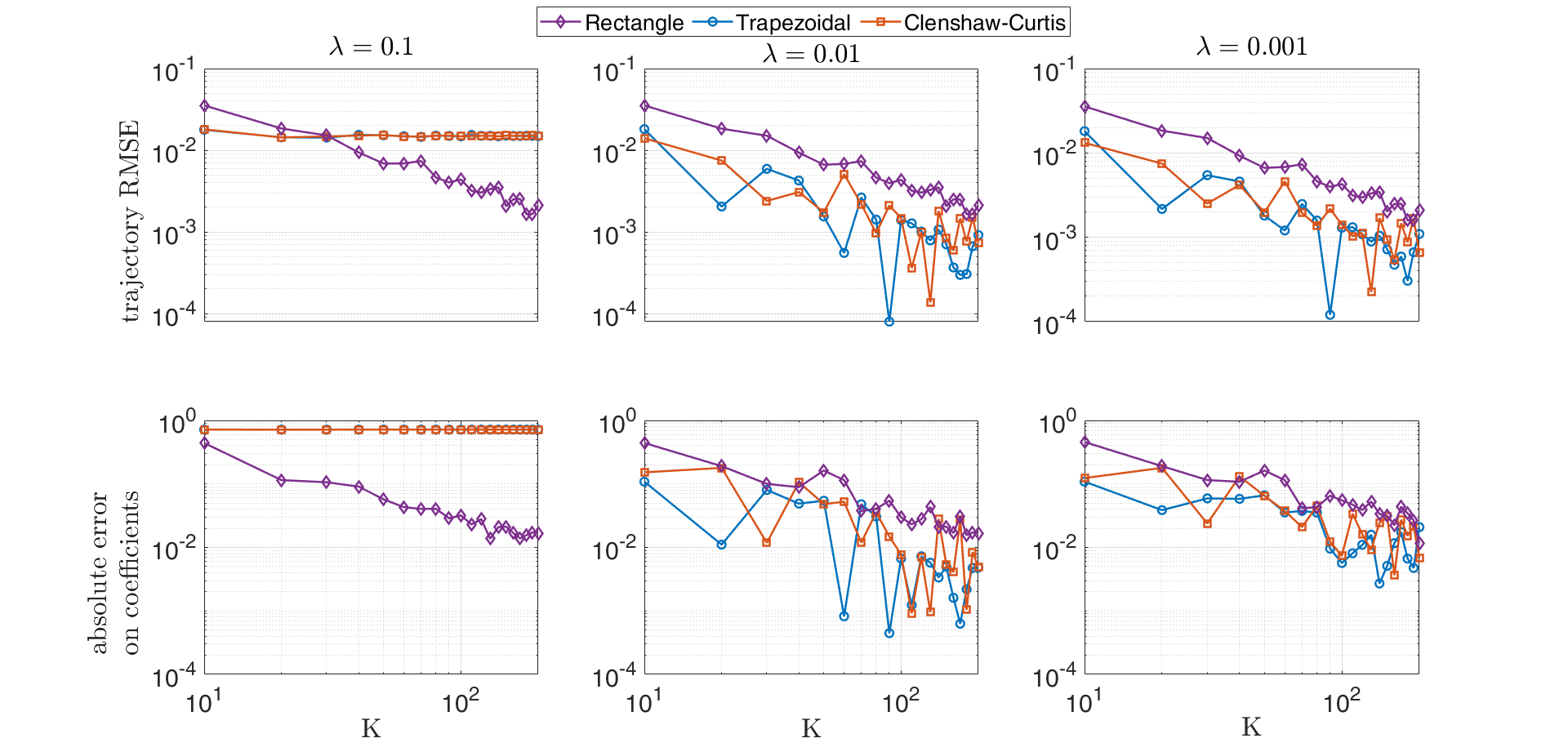}
    \caption{comparison of RMSE on trajectories (top row) and absolute error on coefficients (bottom row) for \eqref{eq:toy_re}, varying $K$ and $\lambda=0.1$ (left), $\lambda=0.01$ (center) and $\lambda=0.001$ (right). $m=50$ uniform samples were collected in the time window $[0,20]$ using $50\%$ for training and $50\%$ for validation.}
    \label{fig:TM2}
\end{figure}

\begin{figure}
    \centering
    \includegraphics[width=1\textwidth]{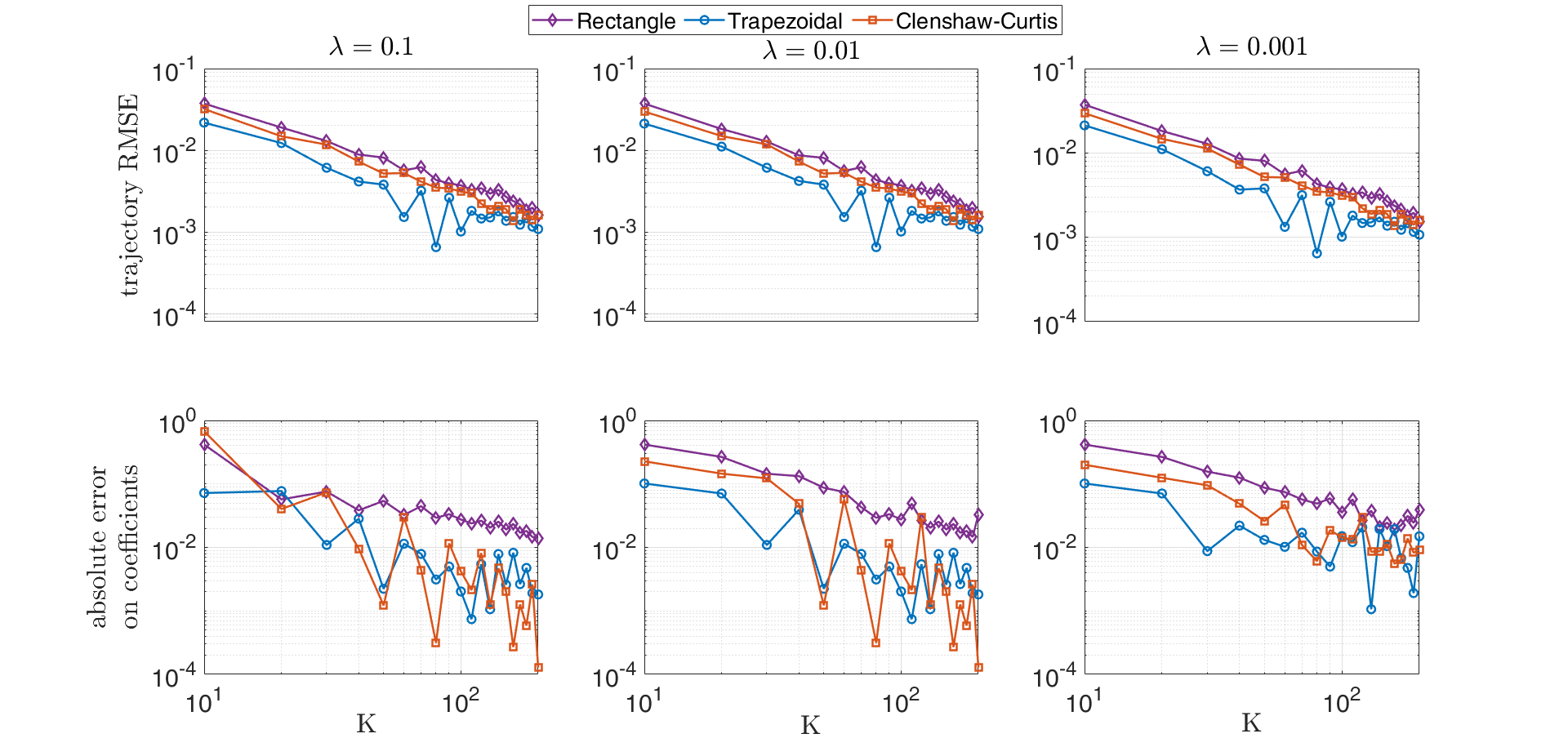}
    \caption{comparison of RMSE on trajectories (top row) and absolute error on coefficients (bottom row), varying $K$ and $\lambda=0.1$ (left), $\lambda=0.01$ (center) and $\lambda=0.001$ (right). $m=200$ uniform samples were collected in the time window $[0,20]$ using $50\%$ for training and $50\%$ for validation.}
    \label{fig:TM200}
\end{figure}

\begin{figure}
    \centering
    \includegraphics[width=1\textwidth]{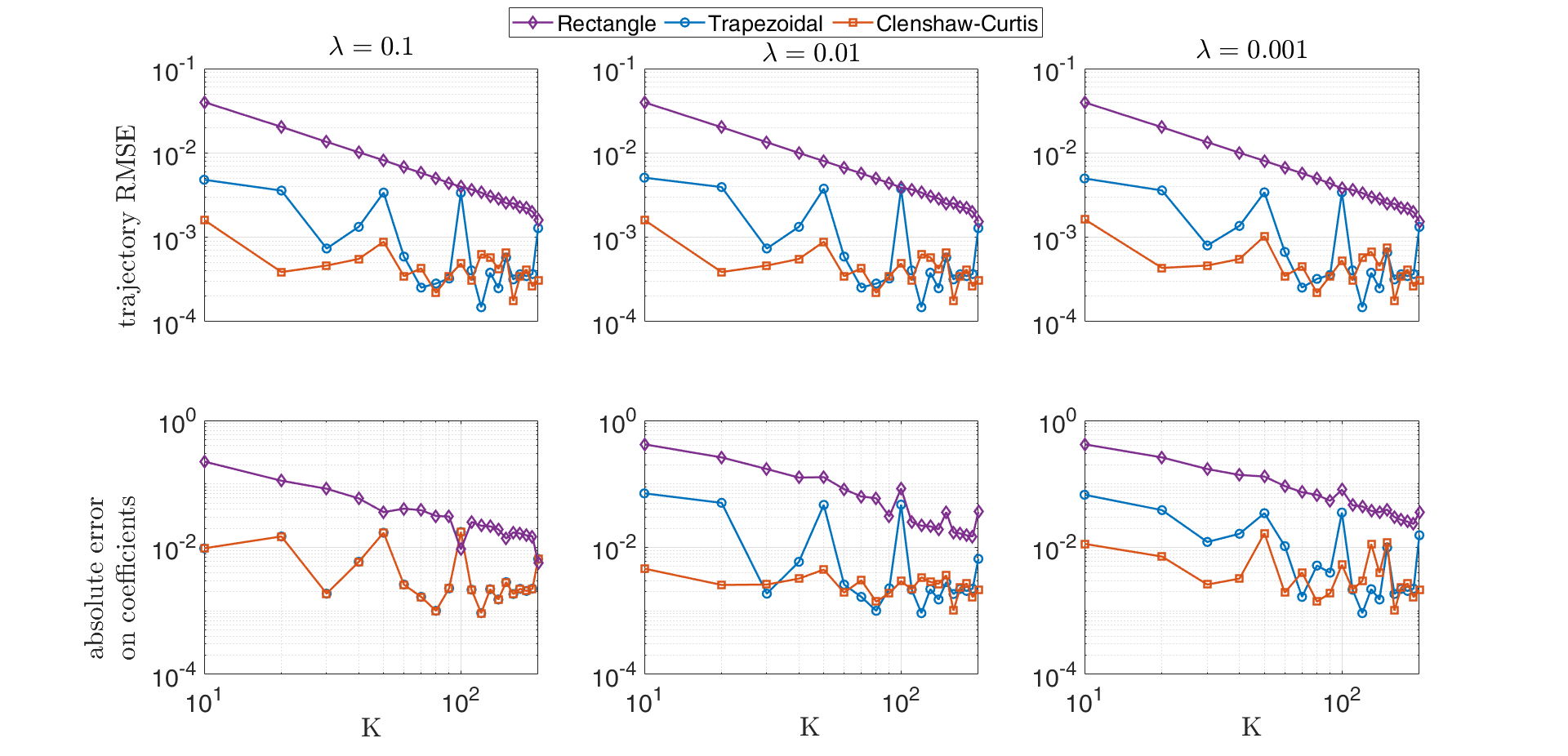}
    \caption{comparison of RMSE on trajectories (top row) and absolute error on coefficients (bottom row), varying $K$ and $\lambda=0.1$ (left), $\lambda=0.01$ (center) and $\lambda=0.001$ (right). $m=1000$ uniform samples were collected in the time window $[0,20]$ using $50\%$ for training and $50\%$ for validation.}
    \label{fig:TM3}
\end{figure}

\begin{figure}
    \centering
    \includegraphics[width=1\textwidth]{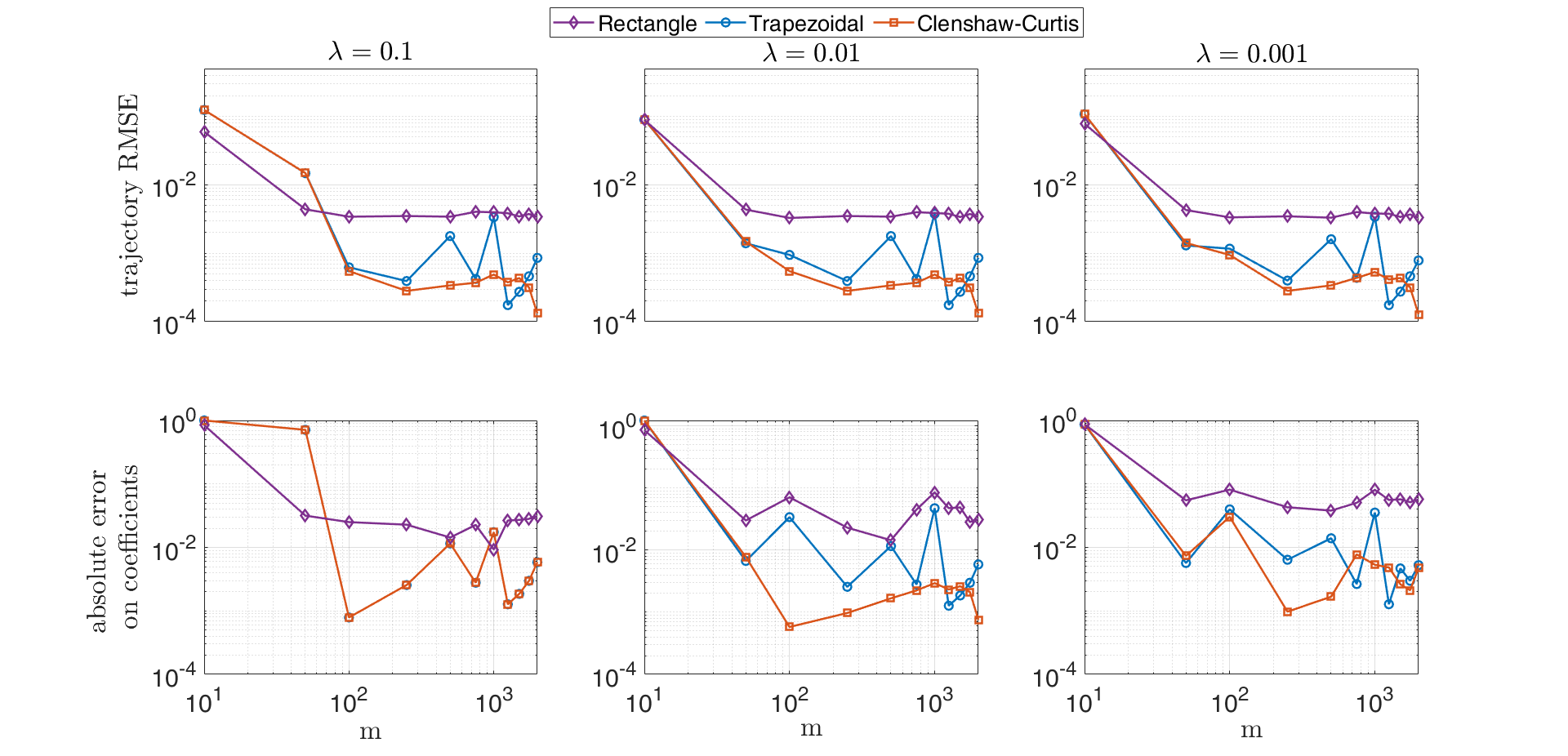}
    \caption{comparison of RMSE on trajectories (top row) and absolute error on coefficients (bottom row), varying $m$ and $\lambda=0.1$ (left), $\lambda=0.01$ (center) and $\lambda=0.001$ (right). $K=100$ quadrature nodes were used and $m$ uniform samples were collected in the time window $[0,20]$ using $50\%$ for training and $50\%$ for validation.}
    \label{fig:TM4}
\end{figure}

Further tests were carried out by varying the number $K$ of quadrature nodes, the hyperparameter $\lambda$ promoting sparsity and the number $m$ of samples, employing also rectangles and Clenshaw–Curtis quadrature \citep{clenshaw1960method,gentleman1972implementing,trefethen2009approximation} for the sake of performance comparison, see \Cref{fig:TM2,fig:TM200,fig:TM3,fig:TM4}. In general and unsurprisingly the accuracy improves by increasing the number $K$ of quadrature nodes independently of the fixed value of $\lambda$ and of the fixed number of samples $m$, see Figure \ref{fig:TM2} for $m=50$,  Figure \ref{fig:TM200} for $m=200$ and Figure \ref{fig:TM3} for $m=1\,000$. The behavior of different quadrature formulas is similar for low or moderate values of $m$, while for larger values differences emerge, with Clenshaw-Curtis prevailing on the other choices. At fixed number of quadrature nodes $K$ the accuracy increases with the number of samples $m$, apparently reaching a barrier already with moderate values on the order of hundreds. Clenshaw-Curtis and trapezoidal rules perform a little better than rectangles, even though no clear trend is visible along the barrier.

The accurate identification of distributed delay kernels via DD-SINDy turns out to be a delicate balance among the adopted quadrature, an appropriate choice of $\lambda$ and an adequate sampling (i.e., $m$). A thorough analysis is out of the scope of the present work, even though the linear decrease of the errors w.r.t. $K$  for all type of quadratures certainly deserves a focused investigation in the future. Let us just notice indeed that while the classical quadrature error is the result of replacing the integrand with a suitable finite-dimensional approximation, typically (piecewise) polynomial on given sets of nodes, in the case of sparse regression one has to take into account also (and mainly) for the linear approximation through the library candidate functions. 

\subsection{An advanced Ricker-type model}\label{sec:adv_ricker}
We consider now a version of \eqref{eq:ricker} more sophisticated than that discussed in Section \ref{sec:ricker}. Indeed, by keeping the same choices of parameters that led therein to \eqref{eq:rickerg1}, let us change $d_{1}$ and $a$ respectively to $d_{1}=0.5$ and $a=\pi/10$, yielding the kernel \eqref{eq:rickerg}
\begin{equation}\label{eq:rickerg2}
g(\sigma,\varphi(\sigma))\coloneqq\gamma\cdot\sigma^{3} e^{4.5\sigma}e^{-\pi\varphi(\sigma)/10}\varphi(\sigma).
\end{equation}
Therefore while the factor $\sigma^{3}$ can be recovered via a polynomial library that includes the candidate function $\sigma$, the two exponential factors need for external optimization of both non-multiplicative parameters $d_{1}$ and $a$ when treating them as unknown. In particular Algorithm \ref{alg:adaptive-distributed-sindy} was used to externally optimize all the unknown values $\tau=1$, $n=4$, $\alpha+d_{1}=4.5$ and $a=\pi/10$. The recovered value for $n$ was rounded to the nearest integer and combined with the optimized $\tau$ to get first $\alpha=n/\tau$ and then $d_{1}$ from the optimized value of $\alpha+d_{1}$, while $a$ was optimized independently.

\begin{figure}[h!]
\centering
\includegraphics[width=0.9\textwidth, height=0.32\linewidth]{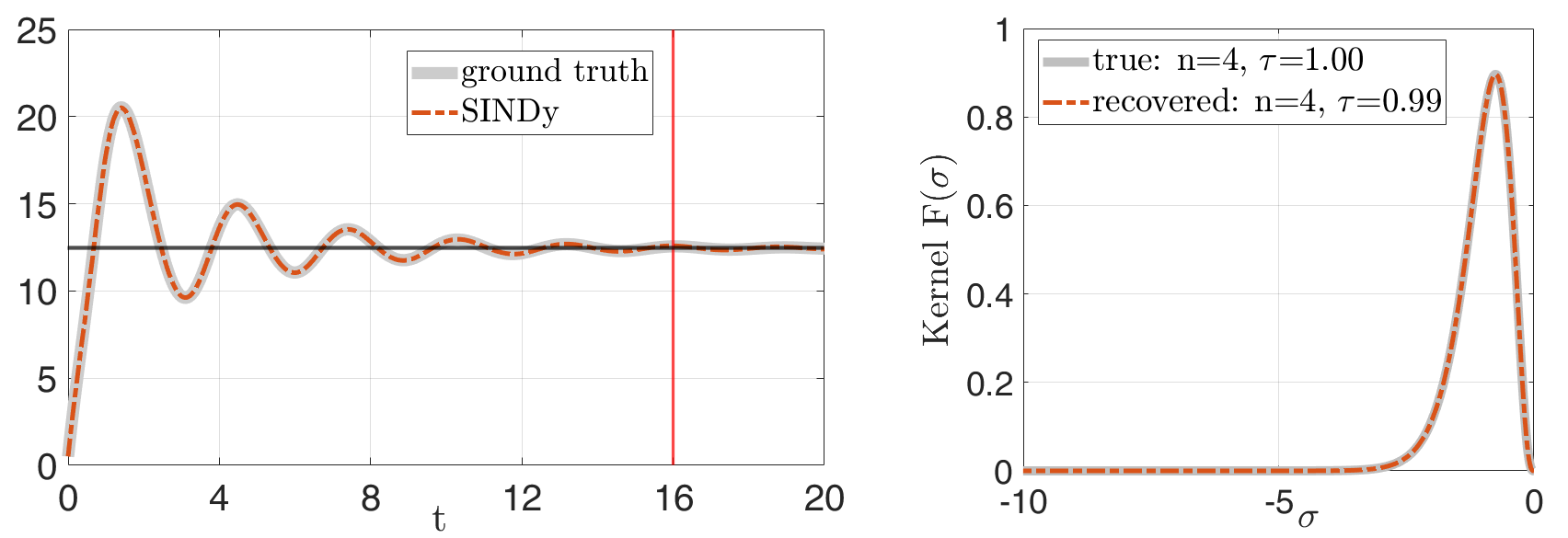}
\caption{comparison of trajectories (left) and kernel functions over the delay interval $[-10\tau,0]$ (right) of \eqref{eq:ricker} with \eqref{eq:rickerg2} by DD-SINDy; corresponding data in Table \ref{tab:comparison_ricker}, horizontal line in the left panel indicates the equilibrium}.\label{fig:ricker_advanced}
\end{figure}

\begin{figure}[h!]
\centering
\includegraphics[width=1\textwidth]{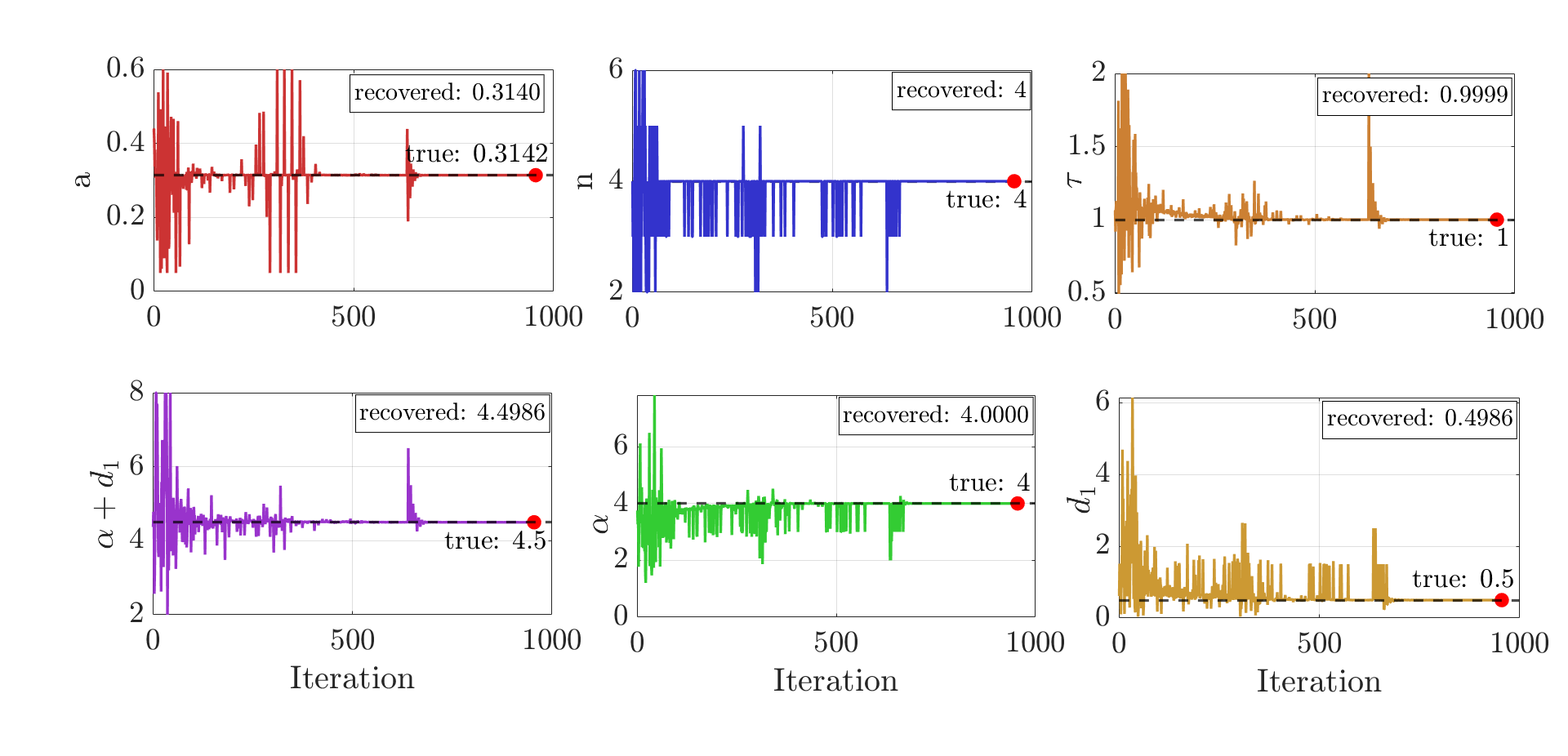}
\caption{iterations of the external optimization of the unknown parameters $\widehat{a}$, $\widehat{n}$, $\widehat{\tau}$, $\widehat{\alpha}+\widehat{d}_1$, $\widehat{\alpha})$ and $\widehat{d}_1$ using particle swarm (PS) for the model \eqref{eq:ricker} with \eqref{eq:rickerg2}; red circles mark final optimized values, corresponding data in Table \ref{tab:comparison_ricker}.}\label{fig:ricker_adv_comp}
\end{figure}

Data for the tests were generated similarly as illustrated in Section \ref{sec:ricker}, collecting now $m=500$ equidistant samples on $[0, 20]$, using the first $80\%$ for training and the rest for validation. A polynomial SINDy library of degree $d=5$ was employed in combination with trapezoidal quadrature with $K=100$ nodes and $\lambda=10^{-2}$. Table \ref{tab:comparison_ricker} collects the results concerning the externally optimized values (top part), the absolute error on the recovered sparse coefficients (central part) and the RMSE on the derivative both for training and validation (bottom part). Correspondingly, Figure \ref{fig:ricker_advanced} shows the recovered trajectory and kernel, while Figure \ref{fig:ricker_adv_comp} illustrates the iterations of the external optimization by Algorithm \ref{alg:adaptive-distributed-sindy} as explained above. 

The obtained results show the excellent performance of DD-SINDy, even when coupled to external optimization to recover unknown delay(s) and non-multiplicative parameters in general. The reconstructed trajectory closely resembles the ground truth during both training and validation. Also the extracted kernel agrees strongly with the actual one. The iterations of Algorithm \ref{alg:adaptive-distributed-sindy} demonstrate effective convergence and alignment with the true values. Optimization took approximately $1.50 \times 10^{3}$ seconds for $971$ internal calls to DD-SINDy to reach a given tolerance of $10^{-4}$, ensuring that PS efficiently identified high-dimensional dynamical structures.
\begin{table}[htbp]
\centering
\small 
\begin{tabular}{lcccc}
\toprule
error on externally&$|n-\widehat{n}|$ & $|d_{1}-\widehat{d}_{1}|$&$|\tau-\widehat{\tau}|$&$|a-\widehat{a}|$ \\ 
optimized values&$0$& $1.43\times10^{-3}$&$1.85\times10^{-5}$&$3.64\times10^{-4}$\\ 
\midrule
error on library&$|d_{0}-\widehat{d}_{0}|$&$|\gamma-\widehat{\gamma}|$&&\\
coefficients &$5.26\times10^{-3}$ & $1.68\times10^{-2}$&&\\
\midrule
RMSE$_{x^\prime}$& training & validation&&\\
& $1.15\times10^{-3}$ & $2.74\times10^{-4}$&&\\
\bottomrule
\end{tabular}%
\caption{absolute errors of optimized values (top), absolute errors on library coefficients (middle) and RMSE on the derivative (bottom) obtained by DD-SINDy on \eqref{eq:ricker} with kernel \eqref{eq:rickerg2}. $m=500$ uniform samples were collected over $[0,20]$ with $K=100$, library degree $d=5$ and $\lambda=10^{-2}$, using $80\%$ for training and $20\%$ for validation.}
 \label{tab:comparison_ricker}
 \end{table}
\footnote{equilibrium $x^*$ satisfies $x'(t)=0$, reducing the RHS of \eqref{eq:ricker} to zero, the horizontal line corresponds to the resulting non-trivial positive solution $x^* > 0$.}

\begin{figure}[h!]
\centering
\includegraphics[width=0.9\textwidth, height=0.32\linewidth]{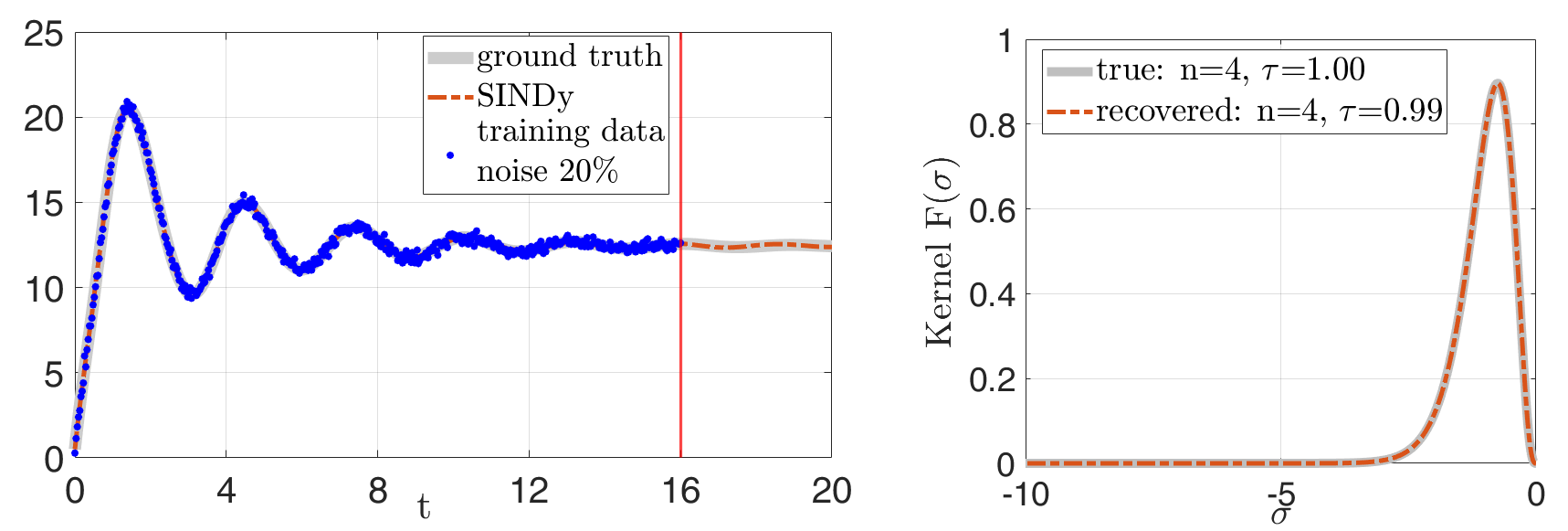}
\caption{comparison of trajectories (left) and kernel functions over the delay interval $[-10\tau,0]$ (right) of \eqref{eq:ricker} with \eqref{eq:rickerg2} by DD-SINDy (noisy data); see Table \ref{tab:comparison_ricker_noise}.}\label{fig:ricker_advanced_noise}
\end{figure}

\begin{figure}[h!]
\centering
\includegraphics[width=1\textwidth]{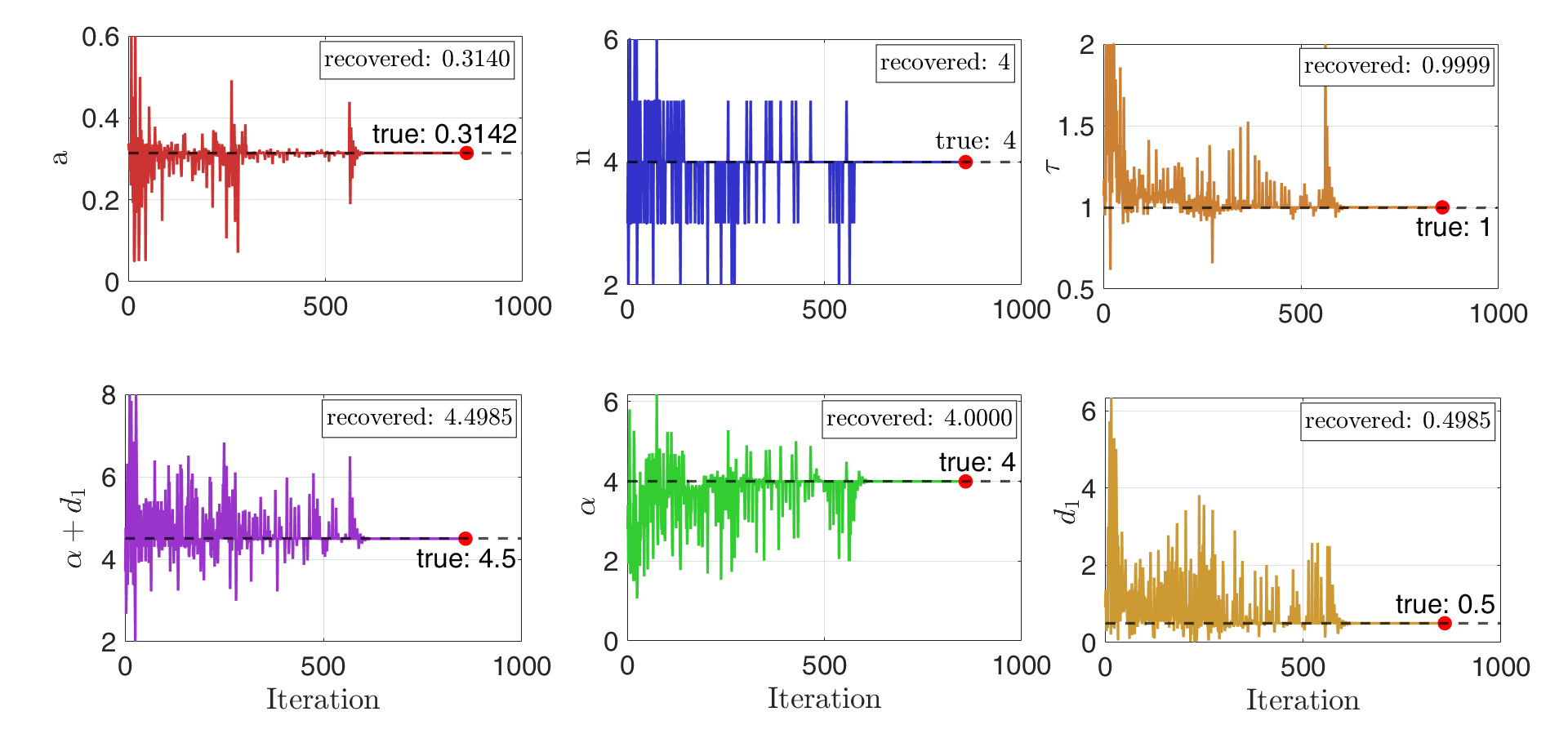}
\caption{iterations of the external optimization of the unknown parameters $\widehat{a}$, $\widehat{n}$, $\widehat{\tau}$, $\widehat{\alpha}+\widehat{d}_1$, $\widehat{\alpha})$ and $\widehat{d}_1$ using particle swarm (PS) for the model \eqref{eq:ricker} with \eqref{eq:rickerg2} on noisy data; red circles mark final optimized values, corresponding data in Table \ref{tab:comparison_ricker_noise}.}\label{fig:ricker_adv_comp_noise}
\end{figure}

The results just described are related to a noise-free dataset. Similar trends were observed by adding $20\%$  white noise to the original dataset, see Table~\ref{tab:comparison_ricker_noise} and Figures \ref{fig:ricker_advanced_noise} and \ref{fig:ricker_adv_comp_noise}. The results based on noisy data demonstrate that DD-SINDy retains high generalization and parameter estimation performance despite noisy conditions.
\begin{table}[htbp]
\centering
\small 
\begin{tabular}{lcccc}
\toprule
error on externally&$|n-\widehat{n}|$ & $|d_{1}-\widehat{d}_{1}|$&$|\tau-\widehat{\tau}|$&$|a-\widehat{a}|$ \\ 
optimized values&$0$& $1.48\times10^{-3}$ & $1.10\times10^{-4}$ & $3.83\times10^{-4}$\\ 
\midrule
error on library&$|d_{0}-\widehat{d}_{0}|$&$|\gamma-\widehat{\gamma}|$&&\\
coefficients& $4.02\times10^{-3}$ & $6.75\times10^{-2}$&&\\
\midrule
RMSE$_{x^\prime}$& training & validation&&\\
&$1.19\times10^{-3}$ & $2.83\times10^{-4}$&&\\
\bottomrule
\end{tabular}
\caption{absolute errors of optimized values (top), absolute errors on library coefficients (middle) and RMSE on the derivative (bottom) obtained by DD-SINDy on \eqref{eq:ricker} with kernel \eqref{eq:rickerg2}. A $20\%$ white noise was added to the $m=500$ uniform samples collected in the time window $[0,20]$ with $K=100$, library degree $d=5$ and $\lambda=10^{-2}$, using $80\%$ for training and $20\%$ for validation.}
 \label{tab:comparison_ricker_noise}
 \end{table}
\subsection{The simplified logistic daphnia model}\label{sec:logistic_daphnia}
As a last case study we consider the coupled RE/DIDE \eqref{eq:daphnia_re_dd}. Fixing $a_{\ast} = 3$, $a_{\dag} = 4$, $r = 1$  and $\gamma = 1$ while varying $\beta$, the consumer free equilibrium $(0,\mathcal{K})$ with $\mathcal{K} = 1$ undergoes a transcritical bifurcation at $\beta=(\mathcal{K}(a_{\dag} - a_{\ast}))^{-1} = 1$. This yields a stable nontrivial equilibrium that undergoes a Hopf bifurcation at $\beta\approx 3.0161$ \citep{breda2020approximation}.

\begin{figure}[h!]
\centering
\includegraphics[width=0.75\linewidth]{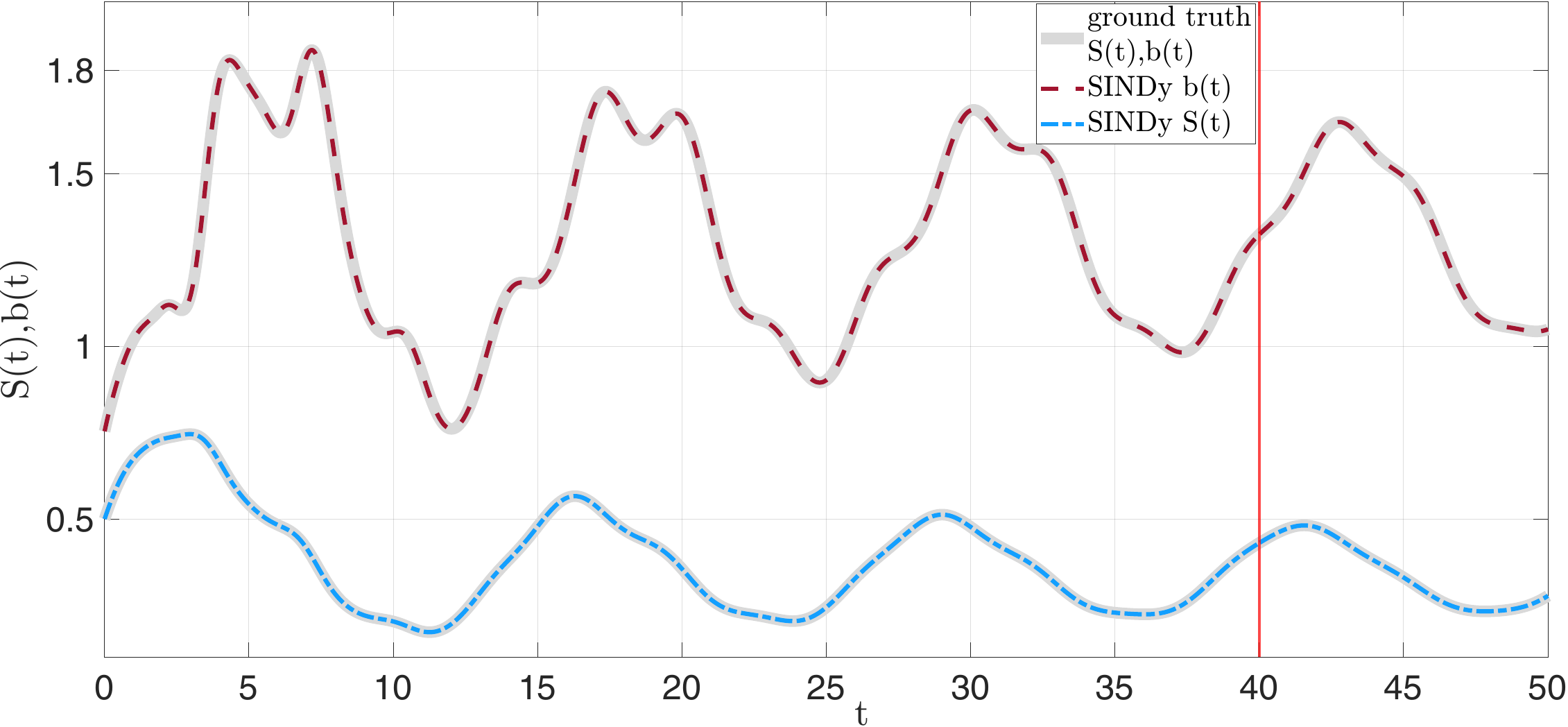}
\caption{comparison of trajectories of \eqref{eq:daphnia_re_dd} reconstructed by DD-SINDy with a polynomial library of degree $d=2$; corresponding data in Table \ref{tab:log_daphnia}.}\label{fig:log_daphnia}
\end{figure}

$m=1\,733$ non-uniform time-step samples for both variables $b$ and $S$ were generated using the MATCONT package \cite{liessi2025new} on the time window $[0,50]$, using $80\%$ for training and $20\%$ for validation. A polynomial SINDy library of degree $d=2$ was employed in combination with trapezoidal quadrature with only $K=8$ nodes and $\lambda=10^{-3}$. Table \ref{tab:log_daphnia} collects the results concerning the externally optimized values (top part), the absolute error on the recovered sparse coefficients (central part) and the RMSE on the derivative both for training and validation (bottom part). Correspondingly, Figure \ref{fig:log_daphnia} shows the recovered trajectories. The integral bounds $a_{\ast}$ and $ a_{\dag}$ were assumed to be unknown, thus identified externally via Algorithm \ref{alg:adaptive-distributed-sindy}.

Again DD-SINDy demonstrates excellent performance also when dealing with coupled REs/DIDEs, even in the case of unknown integration window.

The findings from the logistic Daphnia model \eqref{eq:daphnia_re_dd} underscore two key benefits of the proposed DD-SINDy framework. Firstly, the approach is effective in identifying coupled systems where the dynamics of one state variable (the resource $S$) influence the delayed feedback of another (the consumer birth rate $b$). Secondly, and more importantly, the algorithm adeptly determines the finite integration window $[a_{\ast}, a_{\dag}]$.
The accurate reconstruction of the maturation age $a_{\ast}$ and the maximum age $a_{\dag}$ (refer to Table \ref{tab:log_daphnia}) illustrates that the external optimization loop (Algorithm \ref{alg:adaptive-distributed-sindy}) can identify the 'window of influence' in historical data, providing both a predictive model and biologically meaningful insights into the characteristics of life-history of the population.
\begin{table}[!htbp]
\small
\begin{center}
\begin{tabular}{lcccc}
\toprule
error on externally &{$|a_{\ast}-\widehat{a}_{\ast}|$} & {$|a_{\dag}-\widehat{a}_{\S}|$} & &  \\ 
optimized values& {$1.00\times 10^{-2}$} & {$1.70\times 10^{-2}$}& & \\
\midrule

 & \eqref{eq:daphnia_re_dd_a} & \eqref{eq:daphnia_re_dd_b} & \eqref{eq:daphnia_re_dd_b} & \eqref{eq:daphnia_re_dd_b}  \\
error on library& $S(t)b(t-a)$ & $S(t)$ & $S(t)^{2}$ & $S(t)b(t-a)$ \\
coefficients& $6.05\times 10^{-2}$ & $3.01\times 10^{-3}$ & $4.50\times 10^{-3}$ & $4.44\times 10^{-2}$  \\
 \midrule
RMSE$_{x^\prime}$ & {training} & {validation} & &\\ 
 & {$7.30\times 10^{-4}$} & {$1.40\times 10^{-4}$} & & \\
\bottomrule
\end{tabular}
\end{center}
\caption{absolute errors of optimized values (top), absolute errors on library coefficients (middle) and RMSE on the derivative (bottom) obtained by DD-SINDy on \eqref{eq:daphnia_re_dd}. $m=1\,733$ non-uniform samples were collected in the time window $[0,50]$ with $K=8$, library degree $d=2$ and $\lambda=10^{-3}$, using $80\%$ for training and $20\%$ for validation.}\label{tab:log_daphnia}
\end{table}

\section{Conclusions}\label{sec:conclusion}
This study significantly extends the SINDy framework by expanding its applicability to identify the nonautonomous kernel of the governing equations for systems with distributed delays and REs. By combining a quadrature-based numerical approximation with sparse regression techniques, the proposed method directly reconstructs distributed delay kernels from time-series data. This captures the complex effects of distributed memory inherent in infinite-dimensional dynamical systems.

The proposed methodology introduces several key innovations. First, we develop candidate function libraries that explicitly depend on both delay and state variables, enabling recovery of interpretable kernel structures instead of producing black-box high-dimensional models. Second, the framework applies sparsity-promoting optimization to discover parsimonious representations of distributed memory kernels through weighted sums of quadrature terms. Third, external optimization algorithms jointly identify unknown integration bounds alongside kernel structures, maintaining interpretability throughout. Fourth, comprehensive comparative analysis of various quadrature rules provides practical guidance for practitioners.

This work presents a novel framework marking the first of its kind for distributed delay. Unlike previous research that relied on predefined functional forms to study distributed delays, our method systematically discovers kernel structures and integration limits directly from time-series data. This approach not only enhances interpretability but also surpasses traditional SINDy by revealing complex dynamics within structured population models. This supports future research on more complex functional equations, including systems with features like threshold delays and defined maturation structures \cite{wang2022numerical}. This further broadens the application and scope of data-driven modeling for real-world systems that involve distributed memories.\\




\appendix
\section{Algorithms}\label{sec:appendix}
Pseudocodes for the proposed SINDy framework are given below. Relevant MATLAB demo codes are freely available at \url{https://cdlab.uniud.it/software}.
\begin{algorithm}[H]
\caption{SINDy for Distributed Delay and Renewal equations}
\label{alg:distributed-delay-sindy}
\begin{algorithmic}[1]
\State \textbf{Input:} Data $\mathbf{X},\mathbf{X}'\in\mathbb{R}^{m \times n}$, function library $ \mathcal{F}$, delay interval $[-\tau, 0]$, quadrature rule parameters, regularization parameter $\lambda$
\State \textbf{Output:} Sparse coefficient matrix $\Xi$ representing kernel $g(\sigma,x(t+\sigma))$

\State \textbf{Step 1: Quadrature Setup}
\State Discretize the delay interval $[-\tau, 0]$ into $K$ nodes $\sigma_{k}$ with weights $w_{k}$(e.g., via trapezoidal, Clenshaw-Curtis, etc.)

\State \textbf{Step 2: Data Preparation}
\For{$k = 1$ to $K$}
    \State Obtain time shifted data $\mathbf{X}_{\sigma_{k}}$ by interpolating $\mathbf{X}$ at times $t_{i}+\sigma_{k}$
\EndFor

\State \textbf{Step 3: Library Construction}
   $\Theta = \sum_{k=1}^{K} w_k \cdot \Theta(\sigma_k, \mathbf{X}(t + \sigma_k))$



\State \textbf{Step 4: Sparse Regression}
\For{$j = 1$ to $n$ (each state component)}
    \If{DIDE:} $Y_{j} = \mathbf{X}'_{j}$
    \Else \textbf{ (RE:)} $Y_{j} = \mathbf{X}_{j}$
    \EndIf
    \State Solve sparse regression using STLS or similar:

      \State \quad $\xi_{j} = \arg\min_{\xi} \left(\left\| Y_{j} - \left( \sum_{k=1}^{K} w_{k} \mathbf{\Theta}_k \right) \xi \right\|_{2}+\lambda\|\xi\|_{1} \right) $

\EndFor
\State Assemble coefficient matrix: $\Xi := (\xi_{1}, \xi_{2},\ldots,\xi_n)$

\State \textbf{Step 5: Model Validation}
\State \textbf{return} $\Xi$, validation metrics
\end{algorithmic}
\end{algorithm}

\begin{algorithm}[H]
\caption{SINDy with Parameter Optimization for Distributed Delays}
\label{alg:adaptive-distributed-sindy}
\begin{algorithmic}[1]
\State \textbf{Input:} Data $X$, $\mathbf{X}'$, function library $ \mathcal{F}$, parameter search spaces
\State \textbf{Output:} Optimal parameters $\rho^*$ and sparse coefficient matrix $\Xi^*$

\Function{$\mathcal{J}$}{$\rho$}
    \State Parse parameters: delay bounds $[\tau_{1}, \tau_{2}]$, quadrature settings from $\rho$
    \State Generate nodes $\{\sigma_{k}\}_{k=1}^K$ with weights $\{w_{k}\}_{k=1}^K$ in $[\tau_{1}, \tau_{2}]$ and obtain $\{\mathbf{X}_{\sigma_{k}}\}_{k=1}^K$ by interpolation
    
    \State Construct weighted library matrix $\Theta$ as in Algorithm 1
    
    \For{$j = 1$ to $n$} Solve sparse regression to get $\xi_{j}(\rho)$ 
    \EndFor
    
    \If{DIDE:} $\varepsilon(\rho) = \|\mathbf{X}' - \Theta(\rho)\Xi(\rho)\|_{2}$
    \Else \textbf{(RE:)} $\varepsilon(\rho) = \|X - \Theta(\rho)\Xi(\rho)\|_{2}$
    \EndIf
    \State \textbf{return} $\varepsilon(\rho)$
\EndFunction

\State \textbf{External Parameter Optimization} (PS, BO, etc.): $\rho^* = \arg\min_{\rho} \mathcal{J}(\rho)$ \textbf{return} $\rho^*$, $\Xi^*$
\end{algorithmic}
\end{algorithm}


\section*{Acknowledgements}
DB and MT are members of INdAM research group GNCS; DB is a
member of UMI research group ``Modellistica socio-epidemiologica''. JW is a member of the ``Centre of Excellence in Artificial Intelligence for Public Health Advancement'' and a member of the ``Laboratory for Industrial and Applied Mathematics (LIAM)''. The work of DB was partially supported by the Italian Ministry of University and Research (MUR) through the PRIN 2022 project (No. 20229P2HEA) ``Stochastic numerical modelling for sustainable innovation'', Unit of Udine (CUP G53C24000710006). The work of MT was supported by the Italian Ministry of University and Research (MUR) through a PhD grant PNRR DM351/22 (CUP: G23C22001320003). The work of JW was supported in part by the York Research Chair program (grant number: 492108) and by the NSERC-Sanofi Alliance program in Vaccine Mathematics, Modelling, and Manufacturing (517504).

\end{document}

\endinput